\documentclass[11pt,twoside]{article}
\usepackage{amscd,amssymb,amscd,amsthm,xypic}
\input xy
\xyoption{all}
\newcommand{\PP }{{\mathbb P}}
\newcommand{\QQ }{{\mathbb Q}}
\newcommand{\CC }{{\mathbb C}}
\newcommand{\ZZ }{{\mathbb Z}}

\newtheorem{defn}{Definition}[subsection]
\newtheorem{thm}{Theorem}[subsection]
\newtheorem{prop}{Proposition}[subsection]
\newtheorem{cor}{Corollary}[subsection]
\newtheorem{rem}{Remark}[subsection]

\newcommand{\barefootnote}[1]
{\begingroup
    \renewcommand{\thefootnote}{}
    \footnotetext{#1}
    \renewcommand{\thefootnote}{\arabic{footnote}}
  \endgroup}

\newcommand{\addressemail}[1]
{\begingroup
    \begin{center}
      \vskip-\baselineskip
      #1
    \end{center}
  \endgroup}

\newcommand{\arxurl}[1]
{ \barefootnote{ {\small e-print archive:}
    {\texttt http://lanl.arXiv.org/abs/#1}}}

\renewcommand{\title}[1]
{\begingroup
    \begin{center}
      \vspace{0.4in}
      \bf\huge
      \addtolength{\baselineskip}{5mm}
      #1
    \end{center}
  \endgroup}

\renewcommand{\author}[1] {\begingroup
    \begin{center}
      \vspace{0.4in}
      \bf
      #1
      \vspace{0.2in}
    \end{center}
  \endgroup}

\newcommand{\address}[1]
{\begingroup
    \begin{center}
      #1
    \end{center}
  \endgroup}

\newcounter{gfirstpage}
\newcounter{glastpage}
\newcounter{secondpage}

\begin{document}

\title{ Virtual Class of Zero Loci and Mirror Theorems}

\arxurl{math.AG/0307386}

\author{Artur Elezi}

\address{Department of mathematics and Statistics,\\ American
University}

\addressemail{aelezi@american.edu}

\markboth{\it Virtual class of zero loci and mirror theorems}{\it
Artur Elezi}
\pagestyle{myheadings}
\thispagestyle{empty}

\begin{abstract} \noindent Let $Y$ be the zero loci of a regular section of a
convex vector bundle $E$ over $X$. We provide a new proof of a
conjecture of Cox, Katz and Lee for the virtual class of the
moduli space of genus zero stable maps to $Y$. This in turn yields
the expected relationship between Gromov-Witten theories of $Y$
and $X$ which together with Mirror Theorems allows for the
calculation of enumerative invariants of $Y$ inside of $X$.
\end{abstract}

\section{Introduction}
Let $X$ be a smooth, projective variety over $\CC$. A vector
bundle $$E\rightarrow X$$ is called {\it convex} if
$H^1(f^*(E))=0$ for any morphism $f:\PP^1\rightarrow X$. Let
$Y=Z(s)\subset X$ be the zero locus of a regular section $s$ of a
convex vector bundle $E$ and let $i$ denote the embedding of $Y$
in $X$. It is the relationship between the Gromov-Witten theories
of $Y$ and $X$ that we study here.
\subsection{Virtual class of the zero loci}
Let $\overline M_{0,n}(X,d)$ be the $\QQ$-scheme that represents
coarsely genus zero, $n$-pointed stable maps
$$(C,x_1,x_2,...,x_n,f:C\rightarrow X)$$ of class $d\in H_2(X,\ZZ)$.
Since $E$ is convex the vector spaces $H^0(f^*(E))$ fit into a
$\QQ$-vector bundle $E_d$ on $\overline M_{0,n}(X,d)$. The section
$s$ of $E$ induces a section $\tilde{s}$ of $E_d$ over $\overline
M_{0,0}(X,d)$ via $\tilde{s}((C,f))=s\circ f$. If $i_*(\beta)=d$,
the map $i:Y\hookrightarrow X$ yields an inclusion $i_{\beta}:
\overline M_{0,0}(Y,\beta)\rightarrow \overline M_{0,0}(X,d)$.
Clearly
$$Z(\tilde s)=\coprod_{i_*(\beta)=d}\overline M_{0,0}(Y,\beta).$$
The map $i_*:H_2(Y,\ZZ)\rightarrow (H_2(X,\ZZ)$ is not injective
in general, hence the zero locus $Z(\tilde s)$ may have more then
one connected component. An example is the quadric surface in
$\PP^3$.

The $\QQ$-normal bundle of $Z(\tilde{s})$ in $\overline
M_{0,0}(X,d)$ is $E_d|_{Z(\tilde{s})}$. Let $c_{\text{top}}$
denote the top chern operator. In Section $2.1$ of this paper we
prove the following theorem
\begin{thm}
{\text For any} $d\in H_2(X,\ZZ)$,
$$\sum_{i_*\beta=d}(i_{\beta})_*[\overline
M_{0,0}(Y,\beta)]^{\text{vir}}=c_{\text{top}}(E_d)\cap [\overline
M_{0,0}(X,d)]^{\text{vir}}.$$
\end{thm}

\begin{rem}
Theorem $0.1.1$ was conjectured in \cite{CKL} for any $n$-pointed
moduli stack. But one of the functorial properties of the virtual
class is that
$$\pi_n^*([\overline M_{0,{n-1}}(Y,\beta)]^{\text
virt})=[\overline M_{0,n}(Y,\beta)]^{\text virt}.$$ It is also
true that ${\pi}_{n+1}^*(E_d)=E_d$ (see Lemma $3.4$ in \cite{E1}.)
It suffices then to prove the theorem for $0$-pointed stable maps.
\end{rem}
\begin{rem}
The first proof of this theorem appeared in \cite{KKP}. It is done
in the category of stacks and uses the Behrend-Fantechi
construction of the virtual fundamental class. The version of this
theorem where the virtual classes are those of Li-Tian essentially
follows from Proposition $3.9$ of \cite{LT}. Our approach is a
mixture of the two: we use an extension of the Li-Tian
construction of the virtual fundamental class in the category of
$\QQ$-schemes as well as intersection-theoretic features of a
$\QQ$-scheme coming from its stack nature.
\end{rem}
\subsection{Mirror Theorems} Theorem $0.1.1$ together with the Mirror
Theorem provide a complete answer to the relationship between the
enumerative invariants of $Y$ and those of $X$.

Let $\{T_0=1,T_1,...,T_r,...,T_m\}$ be a basis of $H^*(Y,\QQ)$
such that $$\text{Span}\{T_1,...,T_r\}=H^2(Y,\QQ).$$ Let
$t_1,...,t_r$ be formal variables and $tT=\sum_{i=1}^{r}t_iT_i$.
We denote by $c$ the cotangent line class on $\overline
M_{0,1}(Y,\beta)$, i.e. the first chern class of the line bundle
whose fiber over moduli point $(C,x_1,f)$ is $T^{\vee}_{C,x_1}$.
Let $e$ be the evaluation map on $\overline M_{0,1}(Y,\beta)$. Let
$\hbar$ be a formal parameter. Then the quantum $\mathcal
D$-module structure of the pure quantum cohomology of $Y$  is
determined by the following formal function (see subsection $6$ of
section $1$ in \cite{G1}):
\begin{equation}
J_Y:=e^{\frac{tT}{\hbar}}\cap \left(1+\sum_{0\neq \beta\in
H_2(Y,\ZZ) }q^{\beta}e_*\left(\frac{[\overline
M_{0,1}(Y,\beta)]^{\text{vir}}}{\hbar(\hbar-c)}\right)\right).
\end{equation}
Note that $\displaystyle{e^{\frac{tT}{\hbar}}}$ acts via its power
series expansion (which is finite). For our purposes $J_Y$ will be
viewed as an element of the Novikov completion
$H_*Y[[t_1,...,t_r,\hbar^{-1}]][[q^{\beta}]]$ of the ring
$H_*Y[[t_1,...,t_r,\hbar^{-1}]]$ along the semigroup of rational
curves $\beta$ in $Y$. This generator encodes Gromov-Witten
invariants and the gravitational descendants of $Y$.

For each stable map $(C,x_1,f)\in \overline M_{0,1}(X,d)$, the
sections of $f^*(E)$ that vanish at $x_1$ form a bundle $E'_d$
that fits into an exact sequence:
\begin{equation}
0\rightarrow E'_d\rightarrow E_d\rightarrow e^*(E)\rightarrow 0
\label{eq: spliting}
\end{equation}

The quantum $\mathcal D$-module structure of the $E$-twisted
quantum cohomology of $X$ \cite{E1}\cite{G1} is determined by the
following formal function:
\[
J_E:=e^{\frac{tT'}{\hbar}}c_{\text{top}}(E)\cap
\left(1+\sum_{0\neq d\in H_2(X,\ZZ)
}q^de_*\left(\frac{c_{\text{top}}({E'}_d)\cap [\overline
M_{0,1}(X,d)]^{\text{vir}}}{\hbar(\hbar-c)}\right)\right)
\]
where $tT'$ and $c$ denote similar expressions to those in $J_Y$.
The Mirror Theorem states that for a large class of smooth
varieties $X$, the generator $J_E$ is computable via
hypergeometric series (for proofs of this theorem and its
variations see \cite{B},\cite{E1},\cite{G1},\cite{LLY1}). While
important in itself, this fact is relevant with respect to the
Gromov-Witten theory of $Y$ only if one can show that $J_E$ is
intrinsically related to $Y$. The basic example is when $X$ is a
projective bundle and $E$ is a direct sum of positive line
bundles. It has been shown in \cite{CK} that in this case
$i_*(J_Y)=J_E$. Coupled with the fact that $J_E$ is computable
this allows for the calculation of (at least some of) the
gravitational descendants of $Y$. This is, for example, one way to
compute the enumerative invariants of the quintic threefold. In
section $2.2$ we use Theorem $0.1.1$ to prove the following
generalization: Assume that the map $i_*: H_2(Y)\rightarrow
H_2(X)$ is surjective. Complete a basis $\{T'_1,T'_2,...,T'_r\}$
of $H^2X$ into a basis
$\{T_1=i^*(T'_1),...,T_r=i^*(T'_r),...,T_m\}$ of $H^2Y$. We extend
the map $i_*: H_*Y\rightarrow H_*X$ to a homomorphism of
completions $i_*:
H_*Y[[t_1,...,t_m,\hbar^{-1}]][[q^{\beta}]]\rightarrow
H_*X[[t_1,...,t_r,\hbar^{-1}]][[q^d]]$ via $i_*(t_k)=0$ for $k>r$
and $i_*(q^{\beta})=q^{i_*(\beta)}$.
\begin{thm}
Assume that $i_*:H_2(Y)\rightarrow H_2(X)$ is surjective. Then
$i_*(J_Y)=J_E$.
\end{thm}

{\it Acknowledgements}. The author would like to acknowledge
helpful discussions with Barbara Fantechi, Tom Graber, and Ravi
Vakil. We would also like to thank the referee for providing
corrections and many helpful suggestions.
\section{Background}
\subsection{Stable maps and Gromov-Witten invariants}
Let $g,n$ be non-negative numbers and $d \in H_2(X, \ZZ)$. {\it A
stable map} of genus $g$ with $n$-markings consists of a nodal
curve $C$, an $n$-tuple $(x_1,x_2,...,x_n)$ of smooth points of
$C$ and a map $f:C\rightarrow X$ that has a finite group of
automorphisms. A stable map $(C,x_1,...,x_n,f)$ is said to {\em
represent the curve class $d$} if $f_*[C]=d$. The moduli functor
that parameterizes such stable maps is a proper Deligne-Mumford
stack $\overline {\mathcal M}_{g,n}(X,d)$ and is coarsely
represented by a $\QQ$-scheme denoted by $\overline M_{g,n}(X,d)$
(see section $1.3$ for a discussion on the category of
$\QQ$-schemes). The expected dimension of this moduli stack is
$(\text{dim} X-3)(1-g)+n-K_X\cdot \beta$. Let $x_k$ be one of the
marked points. The evaluation morphism $e_k:\overline
M_{g,n}(X,d)\rightarrow X$ sends a closed point
$(C,x_1,...,x_n,f)$ to $f(x_k)$. The contangent bundle at $x_k$ is
denoted by $\mathcal L_k$. Its fiber over the closed point
$(C,x_1,...,x_n,f)$ is $T^{\vee}_{C,x_k}$. The forgetful morphism
$\pi_k:\overline M_{g,n+1}(X,d)\rightarrow \overline M_{g,n}(X,d)$
forgets the $k$-th marking and stabilizes the source curve. The
universal stable map over $\overline M_{g,n}(X,d)$ is
\[ \begin{CD}
\overline M_{g,n+1}(X,d)@>e_{n+1}>> X \\
@VV \pi_{n+1} V \\
\overline M_{g,n}(X,d).
\end{CD} \]
The bundle $E_d$ from the introduction may be precisely defined as
$E_d:={\pi_{n+1}}_*e_{n+1}^*(E).$
\subsection{The category of $\QQ$-schemes}
Stable maps have nontrivial automorphisms hence the moduli functor
of stable maps is only locally representable, that is the
universal family exists only \'{e}tale locally. It follows that,
in dealing with virtual fundamental class and Gromov-Witten
invariants, one should work in a category that ``remembers'' the
automorphisms. Stacks and, as we will see shortly, $\QQ$-schemes
are two obvious choices.

The key to proving Theorem $0.1.1$ is the functoriality of the
virtual fundamental class as stated in Proposition $1.3.1$. The
construction of the virtual class in the category of stacks has
been done in \cite{BF}. Due to the lack of a good intersection
theory for Artin stacks at the time, the authors of that paper
were forced to impose a technical hypothesis to carry out that
construction. As a result, they were able to prove a stack version
of the Proposition $1.3.1$ that had a limited scope of
applicability (see Proposition $5.10$ in \cite{BF}). For example,
that version could not be used to study the problem we consider
here. Since this technical hypothesis was later removed in
\cite{K}, one obvious approach was to try to prove the
functoriality of the Behrend-Fantechi virtual class construction
in more generality following the standard framework of the
functoriality of the Gysin map (see section 6.5 in \cite{F}).
However, we found the technicalities of the normal cone
construction and the deformation to the normal cone (see section
5.1 in \cite{F}) in the category of stacks hard to overcome. This
functoriality and the subsequent proof of the Cox-Katz-Lee
conjecture in the category of stacks were completed later in
\cite{KKP}.

The Li-Tian (LT) virtual class construction of \cite{K} is free of
any restrictions. Because of this, the authors of that paper were
able to prove the functoriality of their construction in a more
general form, suitable for our problem.

We found the category of $\QQ$-schemes to be the perfect setting
where one can combine the advantages of the LT virtual class
construction and the well established intersection-theoretic
constructions of stacks. The precise definition of $\QQ$-schemes
is due to Lian-Tian (see definition $5.3$ in \cite{LT}). It is a
straightforward generalization of the notion of $\QQ$-varieties
and $\QQ$-stacks introduced by Mumford in \cite{Mu}.

\begin{defn}
A $\QQ$-scheme is a scheme V together with the following data:
\begin{itemize}
\item A finite collection $(V_{\alpha}, G_{\alpha}, q_{\alpha})$
where $V_{\alpha}$ is a quasi-projective scheme, $G_{\alpha}$ is a
finite group acting faithfully on $V_{\alpha}$ and
$q_{\alpha}:V_{\alpha}/G_{\alpha}\rightarrow V$ is an \'{e}tale
map such that $V=\cup Im (q_{\alpha})$.
\item For each pair of indices $(\alpha,\beta)$ there is similarly
$$(V_{(\alpha,\beta)}, G_{(\alpha,\beta)}=G_{\alpha}\times
G_{\beta},~q_{(\alpha,\beta)}:V_{(\alpha,\beta)}\rightarrow V)$$
together with with equivariant finite \'{e}tale maps
$$p_{\alpha}:V_{(\alpha,\beta)}\rightarrow
V_{\alpha},~p_{\beta}:V_{(\alpha,\beta)}\rightarrow V_{\beta}$$
such that $Im q_{(\alpha,\beta)}=Im q_{\alpha}\cap Im q_{\beta}$
and the map $ q_{(\alpha,\beta)}$ factors through both
$q_{\alpha}$ and $q_{\beta}$ via the above maps.
\item For any triple $(\alpha,\beta,\gamma)$ there exists
$(V_{(\alpha,\beta,\gamma)},G_{(\alpha,\beta,\gamma)},
q_{(\alpha,\beta,\gamma)})$ such that
$G_{(\alpha,\beta,\gamma)}=G_{\alpha}\times G_{\beta}\times
G_{\gamma}$, together with equivariant, finite, \'{e}tale maps
from $V_{(\alpha,\beta,\gamma)}$ to
$V_{(\alpha,\beta)},V_{(\alpha,\gamma)},V_{(\beta,\gamma)}$ which
commute with the maps introduced in the second condition and such
that $$Im q_{(\alpha,\beta,\gamma)}=Im q_{\alpha}\cap Im
q_{\beta}\cap Im q_{\gamma}.$$
\end{itemize}
\end{defn}

The motivation for using such a category comes from these
considerations:

1. Any DM stack is generically a quotient of a scheme by a finite
group (see for example Thm. (6.1) of \cite{LM}).

2. The LT construction of the virtual fundamental class can be
done in the category of $\QQ$-schemes.

If the $\QQ$-scheme $V$ represents a moduli functor, then a point
$x\in V_{\alpha,\beta}$ should be thought of as an automorphism
between objects corresponding to $p_{\alpha}(x)$ and
$p_{\beta}(x)$. In fact it is easy to see that the data of a
$\QQ$-scheme determines a stack. Namely, let $R=\coprod
V_{\alpha,\beta}$ and $U=\coprod V_{\alpha}$. The morphisms
$p_{\alpha}~,p_{\beta}$ induce two \'{e}tale morphisms
$p_1,p_2:R\rightarrow U$ and we get an \'{e}tale groupoid scheme
\[
\xymatrix {R\ar@<0.5ex>[r]^{p_1} \ar@<-0.5ex>[r]_{p_2}&U}
\]
and in turn a Deligne-Mumford stack with atlas $U$ (see the
Appendix of \cite{Vi}). It is obvious that this stack has generic
trivial stabilizers. The underlying space of $V$ does not
determine its $\QQ$-scheme structure. For example, in the case of
orbifolds (i.e. $V_{\alpha}$ is smooth for all $\alpha$), one has
to rule out complex reflections.

The definitions of $\QQ$-sheaves and $\QQ$-complexes on a
$\QQ$-scheme $V$ follow naturally. A $\QQ$-sheaf is the collection
of $G_{\alpha}$-equivariant sheaves $\mathcal G_{\alpha}$ on
$V_{\alpha}$ together with isomorphisms $\mathcal
G_{\alpha}\otimes_{\mathcal O_{V_{\alpha}}}\mathcal
O_{V_{\alpha\alpha'}}\simeq \mathcal G_{\alpha'}\otimes_{\mathcal
O_{V_{\alpha'}}}\mathcal O_{V_{\alpha\alpha'}}$ that satisfy the
usual cocycle condition in triple intersections. The
intersection-theoretic machinery that is available for stacks may
be used for $\QQ$-schemes. Of particular importance to us is the
localized top chern class of a section $s$ of a $\QQ$-vector
bundle $E$ on $V$ (see section $14.1$ of \cite {F}). Its existence
and construction in the category of $\QQ$-schemes follows
routinely from the similar construction in the category of stacks
(see for example (ix) of \cite{K}). The local description is
easily obtained by unwinding the definitions. If $s$ is given
locally by $s_{\alpha}$ then the zero locus $Z(s)$ is a closed
$\QQ$-subscheme of $V$ with charts $Z(s_{\alpha})\cap V_{\alpha}$.
Let $F$ be a pure $n$-dimensional subscheme of $V$. The group
$G_{\alpha}$ acts on the normal cone to $Z(s_{\alpha})\cap
q_{\alpha}^{-1}(F)$ in $q_{\alpha}^{-1}(F)$. The quotients can be
patched together to a $\QQ$-cone $C_{Z(s)/F}$ inside of the
restriction to $Z(s)$ of the $\QQ$-bundle $E$ . Let $i$ be the
zero section of this cone. Then the action of the localized top
chern class of $(E,s)$ on $F$ is $i^![C_{Z(s)/F}]$.

\subsection{The virtual fundamental class and the associativity of
the refined Gysin maps} The moduli spaces of stable maps may
behave badly in families and they may have components whose
dimension is bigger than the expected dimension. There is,
however, a cycle of the expected dimension which is deformation
invariant. It is this cycle which is used as the true fundamental
class for intersection theory purposes. In this section we review
the Li-Tian construction of the virtual fundamental class and a
key lemma about the associativity of the refined Gysin maps.

The virtual fundamental class of a moduli functor is constructed
using solely a choice of a tangent-obstruction complex. Our
interest here is the moduli functor $\mathcal F^d_X$ of 0-pointed,
genus zero, degree $d$ stable maps to $X$. We describe the natural
tangent-obstruction complex of $\mathcal F^d_X$. Let $\eta\in
\mathcal F^d_X(S)$ be represented by the following diagram
\[ \begin{CD}
\mathcal X @>f>> X \\
@VV \pi V \\
S
\end{CD} \]
The deformations and obstructions of $\eta$ are described
respectively by the global sections of the sheaves $\mathcal
T^1\mathcal F^d_X(\eta) :=\mathcal E xt^1_{\mathcal
X/S}([f^*(\Omega_X)\rightarrow\Omega_{\mathcal X/S}],\mathcal
O_{\mathcal X})$ and $\mathcal T^2\mathcal F^d_X(\eta):=\mathcal E
xt^2_{\mathcal X/S}([f^*(\Omega_X)\rightarrow\Omega_{\mathcal
X/S}],\mathcal O_{\mathcal X})$. The natural tangent
obstruction-complex for this moduli problem is $\mathcal
T^{\bullet}_{\eta}:=[\mathcal T^1\mathcal F_X(\eta)\rightarrow
\mathcal T^2\mathcal F_X(\eta)]$ with the zero arrow. This complex
is {\it perfect} in the sense that locally, there is a 2-term
complex of locally free sheaves $\mathcal
E^{\bullet}_{\eta}:=[\mathcal E_{\eta,1}\rightarrow \mathcal E
_{\eta,2}]$ whose sheaf cohomology yields the tangent-obstruction
complex $\mathcal H^{\bullet}\mathcal E^{\bullet}_{\eta}=\mathcal
T^{\bullet}_{\eta}$.

The {\em virtual fundamental class} of $\mathcal
T^{\bullet}_{\eta}$ is denoted here by $[\overline
M_{g,n}(X,d)]^{\text{vir}}$. It is a Chow class in the Chow group
of the coarse moduli space $\overline M_{g,n}(X,d)$. If $\overline
M_{g,n}(X,d)$ is an orbifold, its virtual fundamental class
corresponds to $1\in H^*(\overline M_{g,n}(X,d))$ under
Poincar\'{e} duality .

The key to the proof of Theorem $0.1.1$ is a lemma about the
associativity of the refined Gysin maps. Let us first formulate it
for representable functors. Consider a fibre diagram
\begin{equation}
\begin{CD}
W_0                      @>\delta_0>>                W \\
@VV\alpha_0 V                                         @VV\alpha V\\
T_0                @>\delta>>                        T\\ \label{eq:
diagram}
\end{CD}
\end{equation}
where $\delta$ is a regular embedding. Let $\mathcal N$ be the
normal bundle of $T_0$ in $T$. Assume that $W$ and $W_0$ admit
perfect tangent obstruction-complexes $\mathcal T^\bullet_W$ and
$\mathcal T^\bullet_{W_0}$. They are said to be {\em compatible}
relative to the fibre diagram (\ref{eq: diagram}) if for each
affine scheme $S$ and for any morphism $\eta:S\rightarrow
W_0\subset W$ there is an exact sequence
\begin{equation}
0\rightarrow \mathcal T^1_{W_0}(\eta)\rightarrow \mathcal
T^1_{W}(\eta)\rightarrow ({\alpha}_0\circ \eta)^*\mathcal N
\rightarrow \mathcal T^2_{W_0}(\eta)\rightarrow \mathcal
T^2_{W}(\eta)\rightarrow 0.
\end{equation}

Assume that this compatibility satisfies a {\em technical
condition}. Namely, there exists a short exact sequence of 2-term
complexes
\begin{equation}
0\rightarrow [0\rightarrow {\alpha}_0^*{\mathcal N}]\rightarrow
\tilde{\mathcal E}^{\bullet}_{\eta} \rightarrow \mathcal
E^{\bullet}_{\eta}\rightarrow 0 \label{eq: compatibility}
\end{equation}
such that
\begin{itemize}
\item Its long exact sequence of cohomologies is precisely the
exact sequence of the compatibility.
\item The cohomologies of $\tilde{\mathcal E}^{\bullet}_{\eta}$ and
$\mathcal E^{\bullet}_{\eta}$ yield the tangent-obstruction
complex of $W_0$.
\end{itemize}

\begin{prop}(Proposition $3.9$ of \cite{LT}) Assume that $\mathcal T^\bullet_W$ and
$\mathcal T^\bullet_{W_0}$ are compatible and the technical
condition (\ref{eq: compatibility}) is satisfied. Then
$$\delta^![W]^{\text{vir}}=[W_0]^{\text{vir}}$$
where the virtual cycles are with respect to $\mathcal
T^\bullet_W$ and $\mathcal T^\bullet_{W_0}$.
\end{prop}

The case of interest for us is the localized top chern class
\cite{F} of $(E,s)$, where $E\rightarrow Z$ is a vector bundle and
$s$ is a section of $E$. We use diagram (\ref{eq: diagram}) with
$T_0=W=Z$ and $T$ the total space of $E$. Let $\delta=s_E$ be the
zero section of $E$ and $\alpha=s$. It follows that $W_0=Z(s)$ is
the zero locus of $s$.
\begin{cor} \label{corollary: locchern} With the assumptions of the previous proposition
\begin{equation}
{\alpha_0}_* [Z(s)]^{\text{vir}}=c_{\text{top}}(E)\cap
[Z]^{\text{vir}}
\end{equation}
\end{cor}

\section{The Proofs}
\subsection{The associativity of the refined Gysin maps in the category
of $\QQ$-schemes.} Let $(V_{\alpha},G_{\alpha})$ be a $\QQ$-scheme
with a perfect tangent-obstruction $\QQ$-complex $\mathcal
T^{\bullet}$ which is the cohomology of the $\QQ$-complex
$\mathcal E^{\bullet}$. It has been pointed out in \cite{LT} that
the virtual fundamental class in the category of $\QQ$-schemes can
be constructed as follows: First, one constructs a local virtual
cone $C^{\mathcal E^{\bullet}_{\alpha}}$ using the
tangent-obstruction complex $\mathcal H^{\bullet}\mathcal
E^{\bullet}_{\alpha}=\mathcal T^{\bullet}_{\alpha}$ on
$V_{\alpha}$. This local virtual cone sits inside the vector
bundle $Spec~Sym^{\bullet}(\mathcal E_{2,\alpha}^*)$. The key here
is that the local virtual cones do not depend on $\mathcal
E^{\bullet}$ but only on $\mathcal T^{\bullet}$ (Lemma $3.2$ of
\cite{LT}). Obviously the restrictions of ${\mathcal
E^{\bullet}_{\alpha}}$ and ${\mathcal E^{\bullet}_{\beta}}$ to
$V_{\alpha,\beta}$ yield the same tangent-obstruction complex
$\mathcal T^{\bullet}_{\alpha,\beta}$. It follows that the pull
backs of $C^{\mathcal E^{\bullet}_{\alpha}}$ and $C^{\mathcal
E^{\bullet}_{\beta}}$ to $V_{\alpha,\beta}$ are equivariantly
isomorphic. The $G_{\alpha}$-quotients patch together into the
virtual $\QQ$-cone $C^{\mathcal E^{\bullet}}$ inside the
$\QQ$-vector bundle $Spec~Sym^{\bullet}(\mathcal E_2)$. Finally
the virtual fundamental class is the pull back via the zero
section of $[C^{\mathcal E^{\bullet}}]$.

It is known that $\overline M_{0,0}(X,d)$ has a $\QQ$-scheme
structure so that the functors Hom$(-,\overline M_{0,0}(X,d))$ and
$\mathcal F^d_X$ are equivalent. Here is a brief sketch. Let
$(C,f)\in \overline M_{0,0}(X,d)$ be a stable map. Choose divisors
$H_1,H_2,...,H_r$ so that $f$ intersects each $H_i$ transversally
at $y_{i1},...,y_{id_i}$ (where $d_i:=d\cdot H_i$) and $(\tilde
C,y_{ij})$ has no automorphisms. Now $\text {Aut}(f)$ acts on
$(\tilde C,f)\in \overline M_{0,dr}(X,d)\cap e_{ij}^*(H_i)$ by
permuting the markings. Choose a quasiprojective $\text
{Aut}(f)$-equivariant neighborhood $$U_f\subset \overline
M_{0,dr}(X,d)\cap e_{ij}^*(H_i)$$ of $(\tilde C,f)$ such that all
the stable maps in $U_f$ have no automorphisms and the map that
forgets the markings does not change the source curve. There is an
action of $\text {Aut}(f)$ on the universal stable map
$\eta_f:=(\mathcal C_f,F_f)$ over $U_f$. The classifying map
$U_f/\text{Aut}(f)\rightarrow \overline M_{0,0}(X,d)$ is
\'{e}tale. The neighborhoods $U_f$ satisfy the conditions of a
$\QQ$-scheme.

Now, let $E\rightarrow X$ be a convex vector bundle. Recall that
the vector spaces $H^0(f^*(E))$ fit into a $\QQ$-vector bundle
$E_d$ on $\overline M_{0,0}(X,d)$. The section $s$ of $E$ induces
a section $\tilde{s}$ of $E_d$ over $\overline M_{0,0}(X,d)$ via
$\tilde{s}((C,f))=s\circ f$. Let $(C,f)\in Z(\tilde{s})$. By
shrinking $U_f$ if necessary, the group Aut$f$ acts on the
restriction of the universal family $(\mathcal C_f,F_f)$ over
$U^s_f:=Z(\tilde{s})\cap U_f$ and the classifying map
$U^s_f/\text{Aut}(f)\rightarrow Z(\tilde{s})$ is again \'{e}tale.
Just as in the construction of the virtual class of a
$\QQ$-scheme, the local Gysin diagrams (\ref{eq: diagram}) patch
to a global Gysin diagram and that, with the definitions of the
section $1.2$, Corollary \ref{corollary: locchern} holds in the
category of $\QQ$-schemes. The technical assumption is the same;
the only difference is that the complexes become $\QQ$-complexes.
\subsection{Virtual class of the zero loci} The proof of theorem
$0.1.1$ uses corollary \ref{corollary: locchern}. We need to check
that the {\it technical condition} is satisfied. Let $\eta$ be a
$0$-pointed, genus zero stable map of class $d$ over an affine
scheme $S$ represented by the following diagram
\[ \begin{CD}
\mathcal X @>f>>Y\subset X \\
@VV \pi V \\
S
\end{CD} \]
The deformations and obstructions of $\eta$ are described
respectively by the global sections of the sheaves $\mathcal
T^1\mathcal F^d_X(\eta) :=\mathcal E xt^1_{\mathcal
X/S}([f^*(\Omega_X)\rightarrow\Omega_{\mathcal X/S}],\mathcal
O_{\mathcal X})$ and $\mathcal T^2\mathcal F^d_X(\eta):=\mathcal E
xt^2_{\mathcal X/S}([f^*(\Omega_X)\rightarrow\Omega_{\mathcal
X/S}],\mathcal O_{\mathcal X})$. Since the normal bundle of $Y$ in
$X$ is $E|_Y$ the conormal exact sequence writes
\begin{equation}
0\rightarrow f^*(E^*|_Y)\rightarrow f^*(\Omega_X|_Y)\rightarrow
f^*(\Omega_Y)\rightarrow 0 \label{eq: conormalbundle}
\end{equation}

Recall from \cite{LT} that there is a short exact sequence of
$\mathcal O_{\mathcal X}$-sheaves
\begin{equation}
0\rightarrow W_2\rightarrow W_1\rightarrow f^*\Omega_X\rightarrow
0  \label{eq: seslitian}
\end{equation}
such that:
\begin{itemize}
\item $\mathcal E xt^i_{\mathcal X
/S}([W_1\rightarrow\Omega_{\mathcal X /S}],\mathcal O_{\mathcal X
})$ and $\mathcal E xt^i([W_2\rightarrow 0],\mathcal O_{\mathcal
X})$ vanish for $i\neq 1$.

\item Both $\mathcal E_{\eta,1}=\mathcal E
xt^1_{\mathcal X /S} ([W_1\rightarrow\Omega_{\mathcal X
/S}],\mathcal O_{\mathcal X })$ and $\mathcal E_{\eta,2}=\mathcal
E xt^1([W_2\rightarrow 0],\mathcal O_{\mathcal X})$ are locally
free.
\item The sheaf cohomology of
the complex $\mathcal E^{\bullet}_{\eta}=[\mathcal
E_{\eta,1}\rightarrow \mathcal E_{\eta,2}]$ is the
tangent-obstruction complex of the stable map $\eta$, i.e. there
is an exact sequence $$0\rightarrow \mathcal T^1\mathcal
F^d_X(\eta)\rightarrow \mathcal E_{\eta,1}\rightarrow {\mathcal
E}_{\eta,2}\rightarrow T^2\mathcal F^d_X(\eta)\rightarrow 0.$$
\end{itemize}
We pull back (\ref{eq: seslitian}) exact sequence via (\ref{eq:
conormalbundle}) and obtain the following diagram

\[
\xymatrix{& 0 \ar[d] & 0 \ar[d] & 0 \ar[d] & \\
0 \ar[r] & W_2 \ar[d] \ar[r] & \mathcal A \ar[d] \ar[r] & f^*(E^*) \ar[d] \ar[r] & 0 \\
0 \ar[r] & W_2 \ar[d] \ar[r] & W_1 \ar[d] \ar[r] & f^*(\Omega_X|_Y) \ar[d] \ar[r] & 0 \\
0 \ar[r] & 0 \ar[d] \ar[r] & f^*(\Omega_Y) \ar[d] \ar[r] & f^*(\Omega_Y) \ar[d] \ar[r] & 0 \\
& 0 & 0 & 0 &}
\]

Let $\tilde{\mathcal E}_{\eta,2}:=\mathcal E xt^1_{\mathcal
X/S}([\mathcal A\rightarrow 0],\mathcal O_{\mathcal X})$. We apply
the long exact sequence for $\mathcal E xt$ to various exact
sequences obtained from the above diagram. The middle vertical
sequence yields a short exact sequence $$0\rightarrow [\mathcal
A\rightarrow 0]\rightarrow [W_1\rightarrow \Omega_{\mathcal X
/S}]\rightarrow [f^*(\Omega_Y)\rightarrow \Omega_{\mathcal X
/S}]\rightarrow 0.$$ Its long exact sequence for $\mathcal E xt$
yields $$0\rightarrow \mathcal T^1\mathcal
F^{\beta}_Y(\eta)\rightarrow \mathcal E_{\eta,1}\rightarrow
\tilde{\mathcal E}_{\eta,2}\rightarrow T^2\mathcal
F^{\beta}_Y(\eta)\rightarrow 0.$$ Let $\tilde{\mathcal
E}^{\bullet}_{\eta}:=[\mathcal E_{\eta,1}\rightarrow
\tilde{\mathcal E}_{\eta,2}].$ The last exact sequence says that
the cohomology of $\tilde{\mathcal E}^{\bullet}_{\eta}$ is the
tangent obstruction complex of $\eta$. Next, we apply the long
exact sequence for $\mathcal E xt$ to the top horizontal row of
the diagram and use the conditions for $W_i$'s. We obtain
$$0\rightarrow \mathcal E xt^1_{\mathcal X/S}([f^*(E^*)\rightarrow
0],\mathcal O_{\mathcal X})\rightarrow \tilde{\mathcal
E}_{\eta,2}\rightarrow$$ $$\mathcal E_{\eta,2}\rightarrow \mathcal
E xt^2_{\mathcal X/S}([f^*(E^*)\rightarrow 0],\mathcal O_{\mathcal
X}).$$ But one easily sees that
$$\mathcal E xt^2_{\mathcal X/S}([f^*(E^*)\rightarrow 0],\mathcal
O_{\mathcal X})\simeq \mathcal E xt^1(f^*(E^*),\mathcal
O_{\mathcal X})=0$$ and
$$\mathcal E xt^1_{\mathcal X/S}([f^*(E^*)\rightarrow 0],\mathcal
O_{\mathcal X})\simeq \pi_*f^*(E)$$ as $\mathcal O_S$-sheaves. It
follows that there is an exact sequence
\begin{equation}
0\rightarrow [0\rightarrow \pi_*(f^*(E))]\rightarrow
\tilde{\mathcal E}^{\bullet}_{\eta}\rightarrow \mathcal
E^{\bullet}_{\eta}\rightarrow 0.
\end{equation}
Its long exact sequence of sheaf cohomologies is easily seen to be
\begin{equation}
0\rightarrow \mathcal T^1\mathcal F^{\beta}_Y(\eta)\rightarrow
\mathcal T^1\mathcal F^d_X(\eta)\rightarrow \pi_*f^*E \rightarrow
\mathcal
T^2\mathcal F^{\beta}_Y(\eta) \\
\rightarrow \mathcal T^2\mathcal F^d_X (\eta)\rightarrow 0.
\end{equation}
The technical condition is satisfied.

\subsection{\bf Mirror Theorems} We recall the setup. Let $i$ denote
the embedding of $Y$ in $X$. Assume that the map $i_*:
H_2(Y)\rightarrow H_2(X)$ is surjective. Complete a basis
$\{T'_1,T'_2,...,T'_r\}$ of $H^2X$ to a basis
$\{T_1:=i^*(T'_1),...,T_r:=i^*(T'_r),T_{r+1},...,T_m\}$ of $H^2Y$.
Let $tT=\sum_{i=1}^mt_iT_i$ and $tT'=\sum_{i=1}^rt_iT'_i$ where
$t_1,...,t_n$ are variables. The map $i_*: H_*Y\rightarrow H_*X$
extends to a homomorphism of completions
$$i_*: H_*Y[[t_1,...,t_m,\hbar^{-1}]][[q^{\beta}]]\rightarrow
H_*X[[t_1,...,t_r,\hbar^{-1}]][[q^d]]$$ via $i_*(t_k)=0$ for $k>r$
and $i_*(q^{\beta})=q^{i_*(\beta)}$.

\begin{thm}
Assume that $i_*:H_2(Y)\rightarrow H_2(X)$ is surjective. Then
$i_*(J_Y)=J_E$.
\end{thm}

\begin{proof} Recall that
\[
\displaystyle{J_Y:=e^{\frac{tT}{\hbar}}\cap
\left(1+\sum_{\beta\neq 0}q^{\beta}e_*\left(\frac{[\overline
M_{0,1}(Y,\beta)]^{\text{vir}}}{\hbar(\hbar-c)}\right)\right)}
\]
By the definition of $i_*$ and the projection formula we obtain
\[
\displaystyle{i_*(J_Y)=e^{\frac{tT'}{\hbar}}\cap
\left(e_*(1)+\sum_{\beta\neq
0}q^{\beta}i_*\left(e_*\left(\frac{[\overline
M_{0,1}(Y,\beta)]^{\text{vir}}}{\hbar(\hbar-c)}\right)\right)\right)}
\]
As we have said before any $\beta\in H_2(Y,\ZZ)$ such that
$i_*(\beta)=d$ induces a morphism $i_{\beta}: \overline
M_{0,0}(Y,\beta)\rightarrow \overline M_{0,0}(X,d)$. Consider the
following commutative diagram:
\[
\begin{CD}
{\overline M}_{0,1}(Y,\beta) @>{i_{\beta}}>> {\overline M}_{0,1}(X,d) \\
@VVeV                                         @VVeV\\
Y                 @>{i}>>              X\\
\end{CD}
\]
The line bundle $\mathcal L_1$ on $\overline M_{0,1}(Y,\beta)$ is
the pullback via $i_{\beta}$ of the bundle on $\overline
M_{0,1}(X,d)$. By the projection formula and Theorem $2.0.1$:
\[
i_*\left(\sum_{i_*(\beta)=d}e_*\left(\frac{[\overline
M_{0,1}(Y,\beta)]^{\text{vir}}}{\hbar(\hbar-c)}
\right)\right)=e_*\left(\frac{{\sum_{i_*\beta=d}}(i_{\beta})_*[\overline
M_{0,0}(Y,\beta)]^{\text{vir}}}{\hbar(\hbar-c)}\right)
\]
\begin{equation}
=e_*\left(\frac{c_{\text{top}}(E_d)\cap [\overline
M_{0,0}(X,d)]^{\text{vir}}}{\hbar(\hbar-c)}\right)
\end{equation}

The exact sequence (\ref{eq: spliting}) implies
\begin{equation}
c_{\text{top}}(E_d)=c_{\text{top}}(E'_d)e^*(c_{\text{top}}(E))
\end{equation}
hence by the projection formula
\[
i_*\left(\sum_{i_*(\beta)=d}e_*\left(\frac{[\overline
M_{0,1}(Y,\beta)]^{\text{vir}}}{\hbar(\hbar-c)} \right)\right)
\]
\[=c_{\text{top}}(E)\cap e_*\left(\frac{c_{\text{top}}(E'_d)\cap
[\overline M_{0,0}(X,d)]^{\text{vir}}}{\hbar(\hbar-c)}\right)
\]
Note also that $i_*(1)=i_*([Y])=c_{\text{top}}(E)\cap [X]$. Now
the theorem follows readily.
\end{proof}


\begin{thebibliography}{9}
\bibitem{BF} K. Behrend and B. Fantechi, {\it The intrinsic normal
cone,} Invent. Math., 128(1):45-88, 1997.
\bibitem{B} A. Bertram, {\it Another way to enumerate rational
curves with torus actions}, Invent. Math. {\bf 142} (2000)
487-512.
\bibitem{CKYZ} T.-M. Chiang, A, Klemm, S.-T. Yau, and E. Zaslow,
{\it Local Mirror Symmetry: Calculations and Interpretations,}
Adv. Theor. Math. Phys. {\bf 3} (1999), 495-565.
\bibitem{CK} D. Cox and S. Katz, {\it Mirror Symmetry and
Algebraic Geometry}, Mathematical Surveys and Monographs {\bf 68},
AMS, Providence, RI, 1999.
\bibitem{CKL} D. Cox, S. Katz, Y-P. Lee, {\it Virtual fundamental
classes of zero loci,} in {\it Advances in algebraic geometry
motivated by physics(Lowell, MA, 2000)}, Contemp. Math. {\bf 276},
Amer. Math. Soc., Providence, RI, 2001, 157-166.
\bibitem{E1} A. Elezi, {\it Mirror symmetry for concavex vector
bundles on projective spaces}, International Journal of
Mathematics and Mathematical Sciences {\bf 3} (2003), 159-197.
\bibitem{F} W. Fulton, {\it Intersection theory}, Springer-Verlag,
New York-Berlin Heidelberg, 1984.
\bibitem{G1} A. Givental, {\it A mirror theorem for toric complete intersections}, in {\it Topological field theory,
primitive forms and related topics (Kyoto, 1996)}, Progr. Math.,
{\bf 160}, Birkh\H{a}user, 1998 141-175.
\bibitem{G2} A. Givental, {\it Equivariant Gromov-Witten invariants}, Int. Math. Res. Notices {\bf 13} (1996), 613-663.
\bibitem{K} A. Kresch, {\it Cycle groups for Artin stacks},
Invent. Math. {\bf 138} (1999), 495-536.
\bibitem{KKP} B. Kim, A. Kresch and T. Pantev, {\it Functoriality in intersection theory and a conjecture of
Cox, Katz, and Lee}, J. Pure Appl. Algebra {\bf 179} (April 2003),
no. 1.
\bibitem{LM} G. Laumon and L. Moret-Bailly, {\it Champs
Alg\'{e}briques}, Springer Verlag, 2000.
\bibitem{LLY1} B. Lian, K. Liu, and S.-T.Yau, {\it Mirror principle I}, Asian J. Math. Vol. 1, no. 4 (1997), 729-763.
\bibitem{LLY2} B. Lian, K. Liu, and S.-T.Yau, {\it Mirror principle: A Survey}, Current Development in Mathematics, 1998,
International Press, Cambridge MA, 1998, 35-65.
\bibitem{LT} J. Li and G. Tian, {\it Virtual moduli cycles and Gromov-Witten invariants of algebraic
varieties,} Jour. AMS {\bf 11} (1998), no. 1, 119-174.
\bibitem{Mu} D. Mumford, {\it Towards an enumerative geometry of the
moduli space of curves}, Arithmetic and Geometry II, Progress in
Mathematics 36 (1983), 271-326.
\bibitem{Vi} A. Vistoli, {\it Intersection theory on algebraic
stacks and on their moduli spaces}, Inventiones Math. 97 (1989),
613-670.
\end{thebibliography}
\end{document}